\documentclass[11pt, reqno]{amsart}

\usepackage{amssymb}
\usepackage{amsmath}
\usepackage{mathrsfs}
\usepackage{amsfonts}
\usepackage{color}
\usepackage{vmargin}
\usepackage{amsthm}
\usepackage{graphicx}

\newtheorem{theorem}{Theorem}[section]
\newtheorem{lemma}[theorem]{Lemma}
\newtheorem{propi}{Proposition}[section]
\theoremstyle{definition}

\theoremstyle{remark}
\newtheorem{remark}[theorem]{Remark}

\numberwithin{equation}{section}

\usepackage{color}

\theoremstyle{definition}

\newcommand\supp{\mathop{\rm supp}}

\newcommand{\comment}[1]{}

 \newcommand{\di}{\rm{diam}}

\newcommand{\f}{\frac}

\begin{document}

\title[A New Approach to $(1,1)-$Weak-type Estimate for $S_{\alpha,\psi}$ and $g^{*}_{\lambda,\psi}$ ]{A New Approach to $(1,1)-$Weak-type Estimate for the Littlewood-Paley Operators $S_{\alpha,\psi}$ and $g^{*}_{\lambda,\psi}$ }

\author{Arash Ghorbanalizadeh}
\address{Department of Mathematics, Institute for Advanced Studies in Basic Sciences (IASBS), Zanjan 45137-66731, Iran}

\email{ghorbanalizadeh@iasbs.ac.ir}

\author{Monire Mikaeili Nia  }
\address{Department of Mathematics, Institute for Advanced Studies in Basic Sciences (IASBS), Zanjan 45137-66731, Iran}
\email{moniremikaeili@iasbs.ac.ir}





\begin{abstract}
In this study, $(1,1)-$weak type boundedness of square function $S_{\alpha,\psi}$  is obtained by using Nazarov-Treil and Volberg technique. Also using this result, the $(1,1)-$ weak type boundedness of $g^{*}_{\lambda,\psi}$ operator is investigated.
\end{abstract}

\maketitle

\section{Introduction}
Harmonic analysis is more concerned with different operators applied in studying the regularity of solutions of PDEs, ranging from different linear integral operators such as Calder\'on-Zygmund to  non-linear such as Littlewood-Paley function, and so on.

The Littlewood-Paley operators, including $g_{\psi},~~ S_{\alpha,\psi}$ and $g^{*}_{\lambda,\psi},$ have essential role in harmonic analysis while the problem of strong and  weak type boundedness of Littlewood-Paley operators are studied by many authors. The interested reader may refer to papers \cite{BuiHormozi,CHI, GM,HYY,Ler, Ler5,shi}. Motivated by work of \cite{cody}, we do solve the $(1,1)-$weak type boundedness of square function without using Calder\'on-Zygmund decomposition.

It is well known that Calder\'on-Zygmund decomposition related to doubling property of measures in any space is utilized to obtain weak boundedness of different operators.

Moreover, C. B. Stockdale in \cite{cody}, solved weak type boundedness of Calder\'on-Zygmund operators using Nazarov-Treil and Volberg technique. In this method, he replaced Calder\'on-Zygmund decomposition with Whitney decomposition and gave more simple proof for this problem. We note that the paper of Nazarov-Treil and Volberg is dedicated to solve some problems in nonhomogeneous space with polynomial growth giving new technique (see \cite{Nazarov}).

In this work, the weak type boundedness of Littlewood-Paley square function is obtained using Whitney decomposition. Applying  this result, we also derive $(1,1)-$weak type boundedness of $g^{*}_{\lambda, \psi}$. This new proof benefits us to prove the weak type boundedness of square function with no need to Calder\'{o}n-Zygmund decomposition by which it is frequently utilized in solving the boundedness of weak type operators. 

We outline this work as follows. In Section $2, $ we provide the needed background material. In section $3,$ some auxiliary lemmas have been brought which have essential roles on our main results. In section $4,$ we obtain $(1,1)-$weak type boundedness of square function. In section $5,$ we give the $(1,1)-$boundedness of $g^{*}_{\lambda, \psi}$ operator using the $(1,1)-$weak type boundedness of square function.

{\bf Notation.}  Throughout this paper, $C$ denotes positive constants whose values may change from line to line while they may be independent of the main involved parameters. We use $A\lesssim_{s} B$ to denote that $A\leq C_{s} B$ for a positive constant $C$ depending on $s$. In case $C$ only depends on $n$ as dimension of $\mathbb{R}^{n}$, we just use $A\lesssim B.$  Moreover, we use $A\sim_{} B$  to mean there exist $C_1$ and $C_2$ such that $C_1 A \le B \le C_2 A$. Similarly, $A\sim_{s} B$ stands for $C_{s} A\le B \le c_{s} A$ for some $C_{s}, \,\, c_{s}$ scalars depending on $s.$

\section{Preliminaries}
In this section we list some necessary concepts and points for later use.

We start this section recalling the definition of a cube $Q$ in $\mathbb{R}^{n}$ with center $c$ and side length $r$ denoted by $Q(c,r).$ Moreover, we use the notation $Q(c,ar)$ to mean a cube with the same center $c$ as $Q$ and side length $a$ times enlarged the side length of $Q(c,r).$


 \textbf{The Whitney decomposition of open sets in $\mathbb{R}^{n}$\cite{grafakous}}.
An arbitrary open set in
$\mathbb{R}^{n}$ can be decomposed as a union of disjoint cubes whose lengths are proportional to their
distance from the boundary of the open set.
	\begin{propi}\cite{grafakous} {}
	Let $\emptyset\neq\varOmega \subsetneq \mathbb{R}^{n}$ be open.
	Then there exists a family  of closed cubes $\{Q_{i}\}_{i}$ such that
	\begin{itemize}
		\item
		\[\varOmega=\cup_{i}Q_{i},\]
		and the $Q_{i}$'s have disjoint interiors.
		\item
		\begin{equation}\label{whitney}
			\sqrt{n}r_{i}\leq dist (Q_{i},\varOmega^{c})\leq 4\sqrt{n} r_{i},
		\end{equation}
		\item
		If the boundaries of two cubes $Q_{j}$ and $Q_{k}$ touch, then
		\[\frac{1}{4}\leq \dfrac{r_{j}}{r_{k}}\leq4.\]
		\item
		For a given $Q_{j}$ there exist at most $12^{n}Q_{k}$'s that touch it.
	\end{itemize}
	Here, $r_{i}$ describes the cube side length of $Q_{i}.$
\end{propi}
Note that the coefficients in \eqref{whitney} can be changed proportionally.

Suppose $\psi(x,y)$ to be a
real-valued locally integrable function defined off  the diagonal
$x =y$ in $\mathbb{R}^{2n}$ for which the following conditions for any $x,y,h \in {\mathbb R}^n$ and $0<\delta <1,~~ 0<\gamma < 1$ are satisfied:
\begin{itemize}
\item Size condition:
\begin{equation}\label{eq-sizecondition}
|\psi(x,y)|\lesssim \frac{1}{(1+|x-y|)^{n}} \frac{1}{(1+|x-y|)^{\delta}},
\end{equation}

\item  Smoothness condition:
\begin{equation}\label{eq-smoothcondition1}
|\psi(x,y)-\psi(x+h,y)|\lesssim \frac{1}{(1+|x-y|)^{n}}\Bigl(\frac{|h|}{1+|x-y|}\Bigr)^{\gamma}
\frac{1}{(1+|x-y|)^{\delta}} ,
\end{equation}
whenever $|h|<\frac{1}{2}|x-y|$, and
\begin{equation}\label{eq-smoothcondition2}
|\psi(x,y)-\psi(x, y+h)|\lesssim \frac{1}{(1+|x-y|)^{n}}
\Bigl(\frac{|h|}{1+|x-y|}\Bigr)^{\gamma} \frac{1}{(1+|x-y|)^{\delta}} ,
\end{equation}
whenever $|h|<\frac{1}{2}|x-y|$.
\end{itemize}

For $t>0,$ we define a linear operator $\psi_t$ by
$$
\psi_t f (x)=\frac{1}{t^{n}}\int_{\mathbb{R}^n}\psi\Big(\frac{x}{t},\frac{y}{t}\Big)
f(y)dy,
$$
for $f\in {\mathcal S}(\mathbb{R}^n)$.
For $\alpha>0$ and $\lambda>1,$ the square functions
 $S_{\alpha,\psi}$  and $g^{*}_{\lambda,\psi}$ associated to $\psi$ is defined by
\begin{equation}\label{eq:220411-3}
S_{\alpha, \psi} f (x):=\Big(\iint_{\Gamma_\alpha(x)}|\psi_t f (y)|^2
\f{dydt}{t^{n+1}}\Big)^{\frac12} \quad (x \in {\mathbb R}^n),
\end{equation}
 \begin{equation}
	g_{\lambda,\psi}f(x)=\left( \iint_{\mathbb{R}_{+}^{1+n}}\left( \dfrac{t}{t+|x-y|}\right) ^{n\lambda}|\psi_{t}f(y)|^{2}\frac{dydt}{t^{n+1}}\right) ^{\frac{1}{2}}~~~(x \in \mathbb{R}^{n}),
\end{equation}
respectively, where $\Gamma_\alpha(x)$ denotes the cone of aperture $\alpha>1$
centered at $x$ as
\begin{equation}\label{eq:220403-1111}
\Gamma_\alpha(x)=\{(y,t)\in {\mathbb R}^{n}\times (0,\infty): |x-y|<\alpha t\},
\end{equation}
and $\mathbb{R}_{+}^{1+n}=(0,\infty)\times \mathbb{R}^{n}.$
Moreover,  suppose that $S_{1,\psi}$ is $L^2$-bounded. Then since $\|S_{\alpha,\psi}f\|_{L^2}
=\alpha^{n/2}\|S_{1,\psi}f\|_{L^2}$ (see \cite[Lemma 2.4]{HYY}),  consequently $S_{\alpha,\psi}$ is $L^2$-bounded for any
$\alpha>0$.

Also, note that using \cite[Lemma 2.1]{Ler5} and \cite[Lemma 3.1]{BuiHormozi},
we have
\begin{equation}\label{eq:211107-1}
\|S_{\alpha,\psi}\|_{L^1 \to L^{1,\infty}}
\lesssim \alpha^{n} \|S_{1,\psi}\|_{L^1 \to L^{1,\infty}}.
\end{equation}
 This fact convinces us to just focus on $(1,1)-$weak type boundedness of  $S_{1,\psi}$.
\section{Auxiliary lemmas}
In this section, we present some lemmas which will be essential to prove the main results. The first lemma is proved in \cite[Lemma 2]{cody} while the proof of Lemma \ref{norm1} is essentially inspired by \cite[Lemma 2.8]{HYY}. Also, the proof of Lemma \ref{o} and Lemma \ref{2} is inspired by \cite[Lemma 4.1 - Lemma 4.2- p. 315]{torch}.

\begin{lemma}\cite{cody}\label{cube}
Let $\mu$ be a doubling measure on $\mathbb{R}^{n}$ such that \[
\mu(Q(x,ar)) \le C_{\mu,a}\mu(Q(x,r)),
\]
for all $x \in \mathbb{R}^{n}, r>0,$ and $a>1.$ If $N$ is a positive integer, then
\[
\mu(\cup^{N}_{j=1}Q(x_{j}, ar_{j}))\le C_{\mu,a}\mu(\cup^{N}_{j=1}Q(x_{j}, r_{j})),
\]
for all $x_{1}, ..., x_{N} \in \mathbb{R}^{n}, r_{1}, ..., r_{N}>0,$ and $a>1.$
\end{lemma}
\begin{lemma}\label{2}
The following relation holds
\begin{align}\label{1}
&\left(
\int_{r}^\infty r^n
\left|
\frac{1}{(t+|x-c|)^n}
\left(\frac{t}{t+|x-c|}\right)^{\delta}
\right|^2\frac{dt}{t^{n+1}}
\right)^{\frac12}\notag
\\
 &\qquad \qquad\qquad\qquad\qquad\qquad\lesssim_{n,\delta}
	\sum_{k=1}^\infty
	2^{-k n/2}
	\frac{1}{(2^k r+|x-c|)^n}
	\left(\frac{2^{k}r}{2^kr+|x-c|}\right)^{\delta}.
\end{align}
\end{lemma}
\begin{proof}\label{sigma}
To prove \eqref{1}, we use the following facts. In fact, by a change of variable and concavity of function $t^{\frac{1}{2}}, \left(\left(a+b \right)^{\frac{1}{2}}\leq a^{\frac{1}{2}}+b^{\frac{1}{2}}\right),$ we see
	\begin{align*}
		&	\left( \int^{\infty}_{r}r^{n}\left|
		\frac{1}{(t+|x-c|)^n}
		\left(\frac{t}{t+|x-c|}\right)^{\delta}
		\right|^2\frac{dt}{t^{n+1}}\right) ^{\frac{1}{2}}\\
		&=	\left( \sum^{\infty}_{k=1}\int^{2^{k}r}_{2^{k-1}r}r^{n}	\left|
		\frac{1}{(t+|x-c|)^n}
		\left(\frac{t}{t+|x-c|}\right)^{\delta}
		\right|^2\frac{dt}{t^{n+1}}\right) ^{\frac{1}{2}}\\
		&=\left( \sum^{\infty}_{k=1}\int^{1}_{\frac{1}{2}}r^{n}	\left|
		\frac{1}{(2^{k}ru+|x-c|)^n}
		\left(\frac{ 2^{k}ru}{2^{k}ru+|x-c|}\right)^{\delta}
		\right|^2\frac{2^{k}rdu}{(2^{k}r)^{n+1}u^{n+1}}\right) ^{\frac{1}{2}}\\
		&=\left( \sum^{\infty}_{k=1}2^{-kn}\int^{1}_{\frac{1}{2}}	\left|
		\frac{1}{(2^{k}ru+|x-c|)^n}
		\left(\frac{ 2^{k}ru}{2^{k}ru+|x-c|}\right)^{\delta}
		\right|^2\frac{du}{u^{n+1}}\right) ^{\frac{1}{2}}\\
		&\lesssim_{n} \left( \sum^{\infty}_{k=1}2^{-kn}\left|
		\frac{1}{(2^{k-1}r+|x-c|)^n}
		\left(\frac{ 2^{k}r}{2^{k}r+|x-c|}\right)^{\delta}
		\right|^2\right) ^{\frac{1}{2}}\\
		&\lesssim_{n} \sum^{\infty}_{k=1}2^{-kn}\left|
		\frac{1}{(2^{k}r+|x-c|)^n}
		\left(\frac{ 2^{k}r}{2^{k}r+|x-c|}\right)^{\delta}
		\right|^.\\
	\end{align*}
	where in the last three lines we have used the following facts
	\begin{itemize}
		\item
		The function $\dfrac{1}{2^{k}ru+|x-c|}$ is decreasing on $[\frac{1}{2},1],$ and its maximum is $\dfrac{1}{2^{k-1}r+|x-c|}.$ Consequently, maximum of  $\left(\dfrac{1}{2^{k}ru+|x-c|}\right)^{n}$ is $\left(\dfrac{1}{2^{k-1}r+|x-c|}\right)^{n}$ on $[\frac{1}{2},1].$
		\item	
		The function $\frac{ 2^{k}ru}{2^{k}ru+|x-c|}$ is increasing on $[\frac{1}{2}, 1]$, so its maximum will occur in $u=1.$ Hence, maximum of $\left( \frac{ 2^{k}ru}{2^{k}ru+|x-c|}\right) ^{\delta}$ on $[\frac{1}{2},1]$ will happen at $\left( \frac{ 2^{k}r}{2^{k}r+|x-c|}\right) ^{\delta}.$
		\item
		Noting the point that
		\[\dfrac{1}{2^{k}r+|x-c|}=\dfrac{2}{2^{k+1}r+2|x-c|}< \dfrac{2}{ 2^{k+1}r+|x-c|}.\]
		consequently,
		\[\dfrac{1}{(2^{k}r+|x-c|)^{n}}=\dfrac{2^{n}}{(2^{k+1}r+2|x-c|)^{n}}< \dfrac{2^{n}}{(2^{k+1}r+|x-c|)^{n}}.\]

	\end{itemize}
\end{proof}

\begin{lemma}\cite{HYY}\label{norm1}
	If $f \in L^{1}(\mathbb{R}^{n})$ is supported on $Q(c,r)$ and $\int_{Q(c,r)}f(y)dy=0$ for some $c \in \mathbb{R}^{n}, r>0.$ Then,
	\begin{equation*}
		\int_{\mathbb{R}^{n}\setminus Q(c,6nr)}|S_{1,\psi}f(x)|dx\lesssim \|f\|_{L^{1}(\mathbb{R}^{n})}.
	\end{equation*}
\end{lemma}
\begin{proof}
A change of variable on definition of $S_{1,\psi}$ makes us to conclude that
\[S_{1, \psi }(f)(x)=\left(\int_{0}^{\infty}\int_{|z-x|<t}| \psi_{t}(f)(z)|^{2}\frac{dzdt}{t^{n+1}}\right)^{\frac{1}{2}}=\left(\int_{0}^{\infty}\int_{|z|<t}| \psi_{t}(f)(z+x)|^{2}\frac{dzdt}{t^{n+1}}\right)^{\frac{1}{2}}.\]
Moreover, since $f$ is supported on $Q(c,r)$, it follows that
\[
\psi_{t}(f)(x)=\frac{1}{t^{n}}\int_{\mathbb{R}^n} \psi(\frac{x}{t}, \frac{y}{t})f(y)dy=\frac{1}{t^{n}}\int_{Q(c,r)} \psi(\frac{x}{t}, \frac{y}{t})f(y)dy.
\]	
As $y\in Q(c,r),$ then one can conclude $|y-c|<\frac{\sqrt{n}}{2}r$. Moreover, the assumption $x\in \mathbb{R}^{n}\setminus Q(c,6nr)$ implies to have  $|x-c|>3nr.$ The above mentioned assumption  together with  $|z|\le t,$ it follows that
\begin{equation}\label{eq:211112-1}
	2\left( t+ |x+z-y|\right)> t+ |x-y|+t-|z|> t+ |x-y|,
\end{equation}
and that
\begin{align}\label{eq:211112-2}
	2(t+|x-y|)
	&> t+2(|x-c|-|c-y|)\\
	\nonumber
	&\ge t+|x-c|+3nr-\sqrt{n}r\\
	\nonumber
	&> t+|x-c|.
\end{align}
Putting estimates \eqref{eq:211112-1} and \eqref{eq:211112-2} together, it concludes that
\begin{equation}\label{all}
4\left( t+ |x+z-y|\right)>t+|x-c|.
\end{equation}
In order to estimate $S_{1,\psi}f,$ we deal with the following possible cases
\begin{itemize}
  \item$   3nr/2>|x+z-c|,$
\item $ 3nr/2\le|x+z-c|.$
\end{itemize}

Let us start with
\[
3nr/2>|x+z-c|.
\]
We have the assumption $y\in Q(c,r)$, $|x-c|>3nr$, as well as $|z|\le t$.
Then by following observation, it follows that
\[
3n r <|x-c|
\leq |x+z-c|+|z| < \frac{3n }{2} r + t,
\]
and consequently, one has $t> \frac{3n }{2} r> r$.
By using
\eqref{eq-sizecondition} and
\eqref{all}, we estimate
\begin{align*}
	| \psi_{t}(f)(z+x)|&=
	\lefteqn{
		\left|
		\frac{1}{t^n}
		\int_{Q(c,r)}
		\psi\left(\frac{x+z}{t},\frac{y}{t}\right)f(y)dy
		\right|
	}\\
	&\lesssim
	\int_{Q(c,r)}
	\frac{1}{(t+|x+z-y|)^n}
\left(\frac{t}{t+|x+z-y|}\right)^{\delta}|f(y)|dy\\
	&\lesssim_{n,\delta}
	\int_{Q(c,r)}
	\frac{1}{(t+|x-c|)^n}
	\left(\frac{t}{t+|x-c|}\right)^{\delta}|f(y)|dy.
\end{align*}


Hence,
by above computation and Minkowski's inequality, the assumption $t>r$ as well as  Lemma \ref{2} it follows that
\begin{align*}
|S_{1,\psi}f(x)|&=
		\left(
		\int_{\Gamma_{1}(0)}
		\underbrace{	\left|
			\frac{1}{t^n}
			\int_{Q(c,r)}\psi\left(\frac{x+z}{t},\frac{y}{t}\right)f(y)dy
			\right|^2}_{| \psi_{t}(f)(z+x)|^{2}}\frac{dzdt}{t^{n+1}}
		\right)^{\frac12}\\
	&\lesssim_{n,\delta}
	\left(
	\int_{\Gamma_{1}(0)}\left(	\int_{Q(c,r)}
	\frac{1}{(t+|x-c|)^n}
	\left(\frac{t}{t+|x-c|}\right)^{\delta}|f(y)|dy\right)^{2}\frac{dzdt}{t^{n+1}}
	\right)^{\frac12}\\
	&\lesssim_{n,\delta}
	\|f\|_{L^1}
	\left(
	\int_{|z|\le t,\, r< t}
	\left|
	\frac{1}{(t+|x-c|)^n}
\left(\frac{t}{t+|x-c|}\right)^{\delta}
	\right|^2\frac{dzdt}{t^{n+1}}
	\right)^{\frac12}
	\\
	&\sim_{n,\delta}
	\|f\|_{L^1}
	\left(
	\int_{r}^\infty r^n
	\left|
	\frac{1}{(t+|x-c|)^n}
\left(\frac{t}{t+|x-c|}\right)^{\delta}
	\right|^2\frac{dt}{t^{n+1}}
	\right)^{\frac12}
	\\
	&\lesssim_{n,\delta}
	\|f\|_{L^1}
	\sum_{k=1}^\infty
	2^{-k n/2}
	\frac{1}{(2^k r+|x-c|)^n}
	\left(\frac{2^{k}r}{2^kr+|x-c|}\right)^{\delta}.
\end{align*}
Thus, the estimate for $x \in {\mathbb R}^n$ satisfying $|x+z-c|<3nr/2$ is valid.
Consequently, we get
\begin{align*}
&\int_{|x-c|>3nr}\left| S_{1,\psi}f(x)\right| dx
\\
&\lesssim_{n,\delta} \|f\|_{L^1} \sum_{k=0}^\infty 2^{-k n/2} \int_{|x-c|>3nr}\frac{1}{(2^k r+|x-c|)^n}
\left(\frac{2^{k}r}{2^k r+|x-c|}\right)^{\delta}dx
\\
&\lesssim_{n,\delta}
\|f\|_{L^1}
\sum_{k=0}^\infty
2^{k\delta-\frac{kn}{2}}
\int_{|x-c|>3nr}\frac{dx}{|x-c|^{n+\delta}}
\\
&\lesssim_{n,\delta} \|f\|_{L^1}.
\end{align*}
Note that above series in case $n=1$ is convergent only if  $\delta < \frac{1}{2}$. For $n\geq2,$ the above series is automatically convergent.

Let us discuss the case $3nr/2 \le |x+z-c|$.
Again similar to the previous case, suppose that  $y\in Q(c,r)$, $|x-c|>3nr$ and $|z|\le t$, then

\begin{align*}
|x+z-y|
&\ge |x+z-c|-|y-c|
\\
&\ge \frac{3n }{2} r-\frac{\sqrt{n} }{2} r
\\
&> n r\ge 2 \sqrt{n}|y-c|\ge 2|y-c|.
\end{align*}

So,
\begin{equation}\label{smooth}
|x+z-y| >2|y-c|.
\end{equation}

Thus, we can use the smoothness condition of $\psi$ viewing \eqref{smooth} is satisfied.   We estimate by using the assumption of $\supp f \subset Q(c,r),$  \eqref{eq-smoothcondition2},  the moment condition on $f,$ as well as (\ref{all})
\begin{align*}
|\psi_{t}(f)(z+x)|&= \lefteqn{\left|\frac{1}{t^n} \int_{Q(c,r)}
		\psi\left(\frac{x+z}{t},\frac{y}{t}\right)f(y)dy \right|}
\\
&= \left| \frac{1}{t^n}
\int_{Q(c,r)} \psi\left(\frac{x+z}{t},\frac{y}{t}\right)f(y)dy- \underbrace{\frac{1}{t^n}
\int_{Q(c,r)}\psi\left(\frac{x+z}{t},\frac{c}{t}\right)f(y)dy}_{\text{since}\int_{Q(c,r)}f(y)dy=0}	\right|
\\
&\lesssim \int_{Q(c,r)}
	\frac{1}{(t+|x+z-y|)^{n}}\left(\frac{|y-c|}{t+|x+z-y|}\right)^{\gamma}
\left(\frac{t}{t+|x+z-y|}\right)^{\delta}|f(y)|dy
\\
&\lesssim_{n,\delta}
\int_{Q(c,r)}
\frac{1}{(t+|x-c|)^{n}}
\left(\frac{4(\frac{\sqrt{n}}{2}r)}{4(t+|x+z-y|)|}\right)^{\gamma}
\left(\frac{t}{t+|x-c|}\right)^{\delta}|f(y)|dy\\
&\lesssim_{n,\delta}
\int_{Q(c,r)}
\frac{1}{(t+|x-c|)^{n}}
\left(\frac{2\sqrt{n}r}{|x-c|}\right)^{\gamma}
\left(\frac{t}{t+|x-c|}\right)^{\delta}|f(y)|dy.
\end{align*}
If we comibine this estimate
with Minkowski's inequality
and  $|y-c| \le \frac{\sqrt{n}}{2}r,$ then we have
\begin{align}\label{S}
&\left| S_{1, \psi}f(x)\right|=
\left(\int_{\Gamma_{1}(0)}\underbrace{\left|\frac{1}{t^n}\int_{Q(c,r)}\psi\left(\frac{x+z}{t},\frac{y}{t}\right)f(y)dy\right|^2}_{=|\psi_{t}(f )(x+z)|^{2}}\frac{dzdt}{t^{n+1}}
\right)^{\frac12}\notag
\\
&\lesssim_{n,\delta}\left(
\int_{\Gamma_{1}(0)}\left| \int_{Q(c,r)}
\frac{1}{(t+|x-c|)^{n}}
\left(\frac{2\sqrt{n}r}{|x-c|}\right)^{\gamma}
\left(\frac{t}{t+|x-c|}\right)^{\delta}|f(y)|dy\right| ^{2}
\frac{dzdt}{t^{n+1}}
\right)^{\frac12}\notag
\\
&\lesssim_{n,\delta}
\|f\|_{L^1}
\left(
\int_{\Gamma_{1}(0)}
\left|
\frac{1}{(t+|x-c|)^{n}}
\left(\frac{2\sqrt{n}r}{|x-c|}\right)^{\gamma}
\left(\frac{t}{t+|x-c|}\right)^{\delta}
\right|^2\frac{dzdt}{t^{n+1}}
\right)^{\frac12}\notag
\\
&		\lesssim_{n,\delta}  \|f\|_{L^1} \left[ \left(  \int^{|x-c|}_{0}\int_{|z|<t}\left|   \dfrac{1}{\left(  t+|x-c|\right)  ^{n}}\left(   \frac{2\sqrt{n}r}{|x-c|}\right)   ^{\gamma}\left(   \dfrac{t}{t+|x-c|}\right)   ^{\delta}\right|   ^{2}\frac{dzdt}{t^{n+1}}\right)   ^{\frac{1}{2}}\right.\notag
\\
&\left. \qquad \qquad +\left(  \int^{\infty}_{|x-c|}\int_{|z|<t}\left|   \dfrac{1}{\left( t+|x-c|\right) ^{n}}\left(  \frac{2\sqrt{n}r}{|x-c|}\right)  ^{\gamma}\left(  \dfrac{4t}{t+|x-c|}\right)  ^{\delta}\right|   ^{2}\frac{dzdt}{t^{n+1}}\right)  ^{\frac{1}{2}}\right] \notag
\end{align}
\begin{align}
=: \|f\|_{L^1}\left[  I_{1}+ I_{2}\right].
\end{align}

To discuss $I_{1}$
\begin{align*}
	I_{1}&:=
	\left(
	\int^{|x-c|}_{0}\int_{|z|<t}
	\left|
	\frac{1}{(t+|x-c|)^{n}}
\left(\frac{2\sqrt{n}r}{|x-c|}\right)^{\gamma}
	\left(\frac{t}{t+|x-c|}\right)^{\delta}
	\right|^2\frac{dzdt}{t^{n+1}}
	\right)^{\frac12}
	\\
	&\lesssim_{n,\delta}
	\frac{1}{|x-c|^{n}}
\left(\frac{2\sqrt{n}r}{|x-c|}\right)^{\gamma}
\underbrace{	\left(
	\int^{|x-c|}_{0}	
	\left(\frac{t}{t+|x-c|}\right)^{2\delta}
	\frac{dt}{t}
	\right)^{\frac12}}_{:=J}
	\\
	&\lesssim_{n,\delta}
	\frac{1}{|x-c|^{n}}
	\left(\frac{2\sqrt{n}r}{|x-c|}\right)^{\gamma},
\end{align*}
where $J$ is finite. (See Remark \eqref{J}).
To discuss $I_{2},$ noting that assumption $t>|x-c|$ along with $|x-c|>3nr$ implies $t>r.$ Moreover, thanks to Archimedean property and Lemma \ref{2}  we conclude that

\begin{align*}
I_{2}&=\left(
\int^{\infty}_{|x-c|}\int_{|z|<t}
	\left|
\frac{1}{(t+|x-c|)^{n}}
\left(\frac{2\sqrt{n}r}{|x-c|}\right)^{\gamma}
\left(\frac{t}{t+|x-c|}\right)^{\delta}
\right|^2\frac{dzdt}{t^{n+1}}
	\right)^{\frac{1}{2}}
\\
&\lesssim
\left(
\int^{\infty}_{r}\int_{|z|<t}
\left|
\frac{1}{(t+|x-c|)^{n}}
\left(\frac{2\sqrt{n}r}{|x-c|}\right)^{\gamma}
\left(\frac{t}{t+|x-c|}\right)^{\delta}
\right|^2\frac{dzdt}{t^{n+1}}
\right)^{\frac12}
\\
&\sim
\left(
\int^{\infty}_{ r}
r^{n}\left|
\frac{1}{(t+|x-c|)^{n}}
\left(\frac{2\sqrt{n}r}{|x-c|}\right)^{\gamma}
\left(\frac{t}{t+|x-c|}\right)^{\delta}
\right|^2\frac{dt}{t^{n+1}}
\right)^{\frac12}\\
&\lesssim_{\gamma}
\frac{1}{|x-c|^{\gamma}}
\left( \int_{r}^{\infty}
r^{n}\left|
\frac{1}{(t+|x-c|)^{n}}
\left( \frac{t}{t+|x-c|}\right) ^{\delta}
\right|^2\frac{dt}{t^{n+1}}
\right)^{\frac{1}{2}}\\
&\lesssim_{\gamma}
\frac{1}{|x-c|^{\gamma}}	\sum_{k=1}^\infty
2^{-k n/2}
\frac{1}{(2^k r+|x-c|)^n}
\left(\frac{2^{k}r}{2^kr+|x-c|}\right)^{\delta}.
\end{align*}

Similarly, the above series in $n=1$ is convergent only if $\delta < \frac{1}{2}$. For $n\geq2,$ the series is naturally convergent.

Thus, the estimate for $x \in {\mathbb R}^n$ satisfying $|x+z-c|\ge 3nr/2$ is valid.

Consequently, based on above computation we see
\[\int_{|x-c|>3nr}|S_{1,\psi}f(x)|dx\lesssim_{n,\delta,\gamma}\|f\|_{L^{1}}\left( \int_{|x-c|>3nr}I_{1}dx+\int_{|x-c|>3nr}I_{2}dx\right) \]
\[\lesssim_{n,\delta,\gamma,r}\|f\|_{L^1}
\left( \int_{|x-c|>3nr}\frac{1}{|x-c|^{n+\gamma}}dx+\int_{|x-c|>3nr}\frac{1}{|x-c|^{n+\gamma+\delta}}dx\right) \]
\[
\lesssim_{n,\delta,\gamma,r} \|f\|_{L^1}	\left( \frac{1}{\gamma}\frac{1}{(3nr)^{\gamma}}+\frac{1}{(\gamma+\delta)\left(\frac{1}{(3nr)^{\gamma+\delta}} \right) }\right)
\]
\[
\lesssim_{n,\delta,\gamma,r} \|f\|_{L^1},\]
and this completes the proof.	
\end{proof}
\begin{remark}\label{J}
 We mention that following integral 
\[
J:=	\left(
\int^{|x-c|}_{0}	
\left(\frac{t}{t+|x-c|}\right)^{2\delta}
\frac{dt}{t}
\right)^{\frac12}
\]
is finite. This fact is readily verified by 
\[
J\leq\left( \int^{|x-c|}_{0}\dfrac{t^{2\delta-1} }{|x-c|^{2\delta}}dt\right) ^{\frac{1}{2}}
\]
\[
=\frac{1}{|x-c|^{\delta}}\left( \int^{|x-c|}_{0}t^{2\delta-1}dt\right) ^{\frac{1}{2}}< \infty.\]
\end{remark}

\begin{lemma}\label{o}
	Let $\mathcal{O}$ be an open set in $\mathbb{R}^{n}_{+}$ and for $\alpha\geq 1$ associate to it
	\[\mathcal{U}=\{x\in \mathbb{R}^{n}: M\chi_{\mathcal{O}}(x)>\frac{1}{2\alpha^{n}}\},\]
	where $M$ stands for Hardy-Littlewood maximal operator defined by
	\begin{equation}\label{maximal}Mf(x)=\sup_{t>0}\frac{1}{|B(y,t)|}\int_{B(y,t)}|f(x)|dx,
	\end{equation}
	for all balls $B(y,t)\subset \mathbb{R}^{n}$ centered at $y$ with radius $t$ containing $x.$  Then if $\Gamma_{\alpha}(\mathbb{R}^{n}\setminus \mathcal{U})=\cup_{x \in \mathbb{R}^{n}\setminus \mathcal{U}}\Gamma_{\alpha}(x),$ we have
	
	\begin{itemize}
		\item[(i)] If $(y,t)\in \Gamma_{\alpha}(\mathbb{R}^{n}\setminus \mathcal{U}),$ then
		\begin{equation}\label{add}
			|B(y,t)|\leq 2|B(y,t)\cap (\mathbb{R}^{n}\setminus \mathcal{O})|;
		\end{equation}
		
		\item[(ii)]
		\begin{equation}\label{gamma}
			\Gamma_{\alpha}(\mathbb{R}^{n}\setminus \mathcal{U})\subseteq \Gamma_{1}(\mathbb{R}^{n}\setminus \mathcal{O});
		\end{equation}
		
		\item[(iii)]
		\begin{equation}\label{gammaball}
			|B(y,\alpha t)\cap (\mathbb{R}^{n}\setminus \mathcal{U})|\leq 2\alpha^{n}|B(y, t)\cap (\mathbb{R}\setminus \mathcal{O})|.
		\end{equation}
	\end{itemize}
\end{lemma}

\begin{proof}
	If $(y,t)\in \Gamma_{\alpha}(\mathbb{R}^{n}\setminus \mathcal{U}),$ then there is $x \in \mathbb{R}^{n}\setminus \mathcal{U}$ such that $|y-x|<\alpha t,$ or $x \in B(y, \alpha t).$ Then $ x \in (\mathbb{R}^{n}\setminus \mathcal{U})\cap B(y,\alpha t).$  Now based on $\alpha\geq 1,$ having property $B(y,t)\subseteq B(y,\alpha t),$ one can write
	\begin{align}\label{II}
		\dfrac{|B(y,t)\cap \mathcal{O}|}{|B(y,t)|}&= \alpha^{n}\dfrac{|B(y,t)\cap \mathcal{O}|}{|B(y,\alpha t)|}\notag
		\\
		&\leq \alpha^{n}\dfrac{|B(y,\alpha t)\cap \mathcal{O}|}{|B(y,\alpha a t)|}\notag
		\\
		&\leq \alpha^{n} M\chi_{\mathcal{O}}(x)
		\\
		&\leq\frac{\alpha^{n}}{2\alpha^{n}}=\frac{1}{2}.\notag
	\end{align}
	Knowing the fact that
	\begin{equation}\label{I}\dfrac{|B(y,t)\cap \mathcal{O}|}{|B(y,t)|}+\dfrac{|B(y,t)\cap (\mathbb{R}^{n}\setminus\mathcal{O})|}{|B(y,t)|}=1,
	\end{equation}
	and
	\eqref{II}, we get
	\begin{align}
		\dfrac{|B(y,t)\cap (\mathbb{R}^{n}\setminus\mathcal{O})|}{|B(y,t)|}&=1-\dfrac{|B(y,t)\cap \mathcal{O}|}{|B(y,t)|}\notag
		\\
		&\geq \frac{1}{2}.
	\end{align}
	Hence,
	\[
	|B(y,t)|\leq 2|B(y,t)\cap(\mathbb{R}^{n}\setminus \mathcal{O})|,
	\]
	which (i) is concluded.
	
	Suppose $(y,t)\in \Gamma_{\alpha}(\mathbb{R}^{n}\setminus \mathcal{U}).$ Item $(i)$ implies that there is
	\begin{equation}\label{ad}
		w \in B(y,t)\cap (\mathbb{R}^{n}\setminus \mathcal{O})\neq \emptyset.
	\end{equation}
	To prove \eqref{ad}, by reductio ad absurdum, suppose $B(y,t)\cap (\mathbb{R}^{n}\setminus \mathcal{O})= \emptyset.$ Consequently, by \eqref{add}, we conclude that  \[|B(y,t)|=0~~~~\quad\forall t,\]
	which is a contradiction.
	Hence,
	\eqref{ad} implies that $(y,t)\in \Gamma_{1}(w) \subseteq \Gamma_{1}(\mathbb{R}^{n}\setminus \mathcal{O}).$ So, (ii) is valid.
	
	Using property of a ball measure and \eqref{add}, one can obtain
	\[
	|B(y,\alpha t)\cap (\mathbb{R}^{n}\setminus \mathcal{U})|\leq |B(y,\alpha t)|=\alpha^{n}|B(y, t)|\leq 2\alpha^{n}|B(y, t)\cap (\mathbb{R}\setminus \mathcal{O})|,
	\]
	which shows (iii) holds.
\end{proof}

\begin{lemma}\label{2}
	Suppose $\mathcal{O}$ is an open set of finite measure and let $\mathcal{U}$ be associated to $\mathcal{O}$ as in Lemma \ref{o}. Then for $\alpha \geq 1$
	\begin{equation}\label{3.14}
		\int_{\mathbb{R}^{n}\setminus \mathcal{U}}S_{\alpha, \psi}f(x)^{2}dx\leq 2\alpha^{n} \int_{\mathbb{R}^{n}\setminus \mathcal{O}}S_{1, \psi}f(x)^{2}dx.
	\end{equation}
\end{lemma}
\begin{proof}
	On one hand, from the definition of Lusin area function, Tonelli's theorem, \eqref{gamma} and \eqref{gammaball} respectively, we will have
	\begin{align}\label{a1}
		&\int_{\mathbb{R}^{n}\setminus \mathcal{U}}S_{\alpha, \psi}f(x)^{2}dx\notag
		\\
		&=\int_{\mathbb{R}^{n}\setminus \mathcal{U}}\int^{\infty}_{0}\int_{|x-y|<\alpha t}|\psi_{t} f(y)|^{2}\frac{dydt}{t^{n+1}}dx\notag
		\\
		&=\int_{\mathbb{R}^{n}\setminus \mathcal{U}}\int^{\infty}_{0}\int_{|x-y|<\alpha t, x \in \mathbb{R}^{n}\setminus \mathcal{U}} |\psi_{t} f(y)|^{2}\frac{dydt}{t^{n+1}}dx\notag
		\\
		&=\int_{\mathbb{R}^{n}\setminus \mathcal{U}}\int^{\infty}_{0} \int_{|x-y|<\alpha t, x \in \mathbb{R}^{n}\setminus \mathcal{U}} 1_{(B(y,\alpha t)\cap (\mathbb{R}^{n}\setminus \mathcal{U}))}(x) |\psi_{t} f(y)|^{2}  \frac{dydt}{t^{n+1}}dx\notag
		\\
		&=\int_{\mathbb{R}^{n}\setminus \mathcal{U}} \iint_{\Gamma_{\alpha}(\mathbb{R}^{n}\setminus \mathcal{U})} 1_{(B(y,\alpha t)\cap (\mathbb{R}^{n}\setminus \mathcal{U}))}(x) |\psi_{t} f(y)|^{2}  \frac{dydt}{t^{n+1}}dx\notag
		\\
		&=\iint_{\Gamma_{\alpha}(\mathbb{R}^{n}\setminus \mathcal{U})}|\psi_{t} f(y)|^{2}\underbrace{\left(\int_{\mathbb{R}^{n}\setminus \mathcal{U}} 1_{(B(y,\alpha t)\cap (\mathbb{R}^{n}\setminus \mathcal{U}))}(x)dx\right) }_{|B(y,\alpha t) \cap (\mathbb{R}^{n}\setminus \mathcal{U}))|} \frac{dydt}{t^{n+1}}\notag
		\\
		&=\iint_{\Gamma_{\alpha}(\mathbb{R}^{n}\setminus \mathcal{U})}|\psi_{t} f(y)|^{2}|B(y,\alpha t)\cap (\mathbb{R}^{n}\setminus \mathcal{U})|\frac{dydt}{t^{n+1}}\notag
		\\
		&\leq\iint_{\Gamma_{1}(\mathbb{R}^{n}\setminus \mathcal{O})}|\psi_{t} f(y)|^{2}|B_(y,\alpha t)\cap (\mathbb{R}^{n}\setminus \mathcal{U})|\frac{dydt}{t^{n+1}}\notag
		\\
		&\leq 2\alpha^{n}  \iint_{\Gamma_{1}(\mathbb{R}^{n}\setminus \mathcal{O})}|\psi_{t} f(y)|^{2} |B(y, t)\cap (\mathbb{R}^{n}\setminus \mathcal{O})|\frac{dydt}{t^{n+1}}.
	\end{align}
	On the other hand, similar to the above argument, from the definition of Lusin area function  one can have
	\begin{align}\label{a2}
		&	\int_{\mathbb{R}^{n}\setminus \mathcal{O}}S_{1, \psi}f(x)^{2}dx=\int_{\mathbb{R}^{n}\setminus \mathcal{O}}\left( \int^{\infty}_{0}\int_{|x-y|<t}|\psi_{t} f (y)|^{2} \frac{dydt}{t^{n+1}}\right) dx\notag\\
		&=\int_{\mathbb{R}^{n}\setminus \mathcal{O}}\int^{\infty}_{0}\int_{|x-y|<t, x \in \mathbb{R}^{n}\setminus\mathcal{O}}|\psi_{t} f (y)|^{2} \frac{dydt}{t^{n+1}}dx\notag\\
		&=\int_{\mathbb{R}^{n}\setminus \mathcal{O}}\iint_{\Gamma_{1}(\mathbb{R}^{n}\setminus \mathcal{O})}1_{B(y,t)\cap (\mathbb{R}^{n}\setminus \mathcal{O})}(x)|\psi_{t}f (y)|^{2} \frac{dydt}{t^{n+1}}dx\notag\\
		&=\iint_{\Gamma_{1}({\mathbb{R}^{n}\setminus \mathcal{O}})}|\psi_{t} f (y)|^{2}\underbrace{\left( 	\int_{\mathbb{R}^{n}\setminus \mathcal{O}}1_{B(y,t)\cap (\mathbb{R}^{n}\setminus \mathcal{O})}(x)dx\right)}_{=|B(y,t)\cap (\mathbb{R}^{n}\setminus \mathcal{O})|}  \frac{dydt}{t^{n+1}}.
	\end{align}
	Comparing \eqref{a1} with \eqref{a2}, we can conclude the desired result, \eqref{3.14}.
\end{proof}

\section{$(1,1)-$ weak type boundedness of square function}
In this section, we plan to show the $(1,1)-$weak type boundedness of square function using Whitney decomposition.
\begin{theorem}
	The square function operator \[S_{1,\psi}f(x)=\left(\int^{\infty}_{0}\int_{|x-z|<t}|\psi_{t}f(z)|^{2}\frac{dzdt}{t^{n+1}}\right)^{\frac{1}{2}}\]
		satisfies
	\[
	\|S_{1,\psi}f\|_{L^{1,\infty}}\lesssim \|f\|_{L^{1}(\mathbb{R}^{n})},
	\]
	for all $f \in L^{1}(\mathbb{R}^{n}).$
\end{theorem}
\begin{proof}
	Let $\rho>0$ be given. We intend to verify the following
	\begin{equation*}
		|\{x\in \mathbb{R}^{n}:|S_{1,\psi}f(x)|>\rho\}|\lesssim \frac{1}{\rho
		}\|f\|_{L^{1}(\mathbb{R}^{n})}.
	\end{equation*}
By density, we suppose $f$ is a non-negative continuous function with compact support. Let
\[
\Omega:=\{x\in \mathbb{R}^{n}: f(x)>\rho\}.
\]
Now, apply a Whitney decomposition to write $\Omega:=\{x\in \mathbb{R}^{n}: f(x)>\rho\}=\cup^{\infty}_{i=1}Q_{i},$ a disjoint union of dyadic cubes where
\[
6\sqrt{n}~~ \di(Q_{i}) \le d\left(Q_{i},\mathbb{R}^{n}\setminus \Omega\right) \le 24 \sqrt{n}~~ \di(Q_{i}).
\]

Put
\[g:=f1_{\mathbb{R}^{n}\setminus \Omega},\qquad b:=f1_{\Omega} \qquad \text{and} \qquad b_{i}:=f1_{Q_{i}}.
\]
Then,
\[
f=g+b=g+\sum^{\infty}_{i=1}b_{i},
\]
where
\begin{itemize}
	\item[(1)]
$\|g\|_{L^{\infty}(\mathbb{R}^{n})}\le \rho$ and $\|g\|_{L^{1}(\mathbb{R}^{n})} \le \|f\|_{L^{1}(\mathbb{R}^{n})},$
	\item[(2)]
	the $b_{i}$ are supported on pairwise disjoint cubes $Q_{i}$ with $\sum^{\infty}_{i=1}|Q_{i}| \le \frac{1}{\rho}\|f\|_{L^{1}(\mathbb{R}^{n})},$  and
	\item[(3)]
	$\|b\|_{L^{1}(\mathbb{R}^{n})}\le \|f\|_{L^{1}(\mathbb{R}^{n})}.$
\end{itemize}
We have
\begin{align*}
|\{x\in \mathbb{R}^{n}: \rho <|S_{1, \psi}f(x)|\}| &\le |\{x \in \mathbb{R}^{n}: |S_{1, \psi}g(x)|\geq\frac{\rho }{2\sqrt{c_{2}}}\}|+|\{x \in \mathbb{R}^{n}: |S_{1, \psi}b(x)|\geq\frac{\rho }{2\sqrt{c_{2}}}\}|\\
\end{align*}
We control the first term, using Chebychev's inequality, the $L^{2}-$boundedness of $S_{1,\psi}$ on $\mathbb{R}^{n}$, and property (1) to estimate
\begin{align*}
	|\{x \in \mathbb{R}^{n}: |S_{1, \psi}g(x)|\gtrsim \frac{\rho }{2}\}| &\lesssim \dfrac{\|S_{1, \psi}g(x)\|_{L^{2}}^{2}}{(\frac{\rho}{2})^{2}}	\\
	&\lesssim \dfrac{\int_{\mathbb{R}^{n}}|S_{1, \psi}g(x)|^{2}dx}{\rho^{2}}\\
	& \lesssim \dfrac{\int_{\mathbb{R}^{n}}|g(x)|^{2}dx}{\rho^{2}}\\
	&\le\dfrac{\int_{\mathbb{R}^{n}}|g(x)|dx}{\rho}\\
	& \le \dfrac{\int_{\mathbb{R}^{n}}|f(x)|dx}{\rho}\\ 	
	&=\frac{\|f\|_{L^{1}(\mathbb{R}^{n})}}{\rho}.
\end{align*}


For positive integers $N,$ set 
\[b^{(N)}:=\sum^{N}_{i=1}b_{i}\]
 To control the second term, it suffices to handle
\[\{x\in\mathbb{R}^{n}:S_{1,\psi}b^{(N)}(x)\geq \frac{\rho }{2\sqrt{C_{2}}}\}\]
 
 uniformly in N.

Let $c_{i}$ denote the center of $Q_{i}$ and let $a_{i}:=\int_{Q_{i}}b_{i}(x)dx.$
Set
$E_{1} := Q(c_{1}, r_{1}),$
where $r_{1} > 0$ is chosen so that $|E_{1}|:=\frac{a_{1}}{\lambda}.$ 
In general, for $i = 2, 3, ..., N,$ set
$E_{i} := Q(c_{i}, r_{i})\setminus \cup^{i-1}_{k=1}E_{k},$
where $r_{i} > 0$ is chosen so that $|E_{i}|=\frac{a_{i}}{\lambda}.$ For $i=1, 2, 3, ..., N,$
Set $E_{i}=Q(c_{i},r_{i})$ where $r_{i}>0$ is chosen so that $|E_{i}|=\frac{a_{i}}{\rho}.$
Define
$$
E:=\cup^{N}_{i=1}E_{i}=\cup^{N}_{i=1}Q(c_{i},r_{i})\qquad E^{*}:=\cup^{N}_{i=1}Q(c_{i},6nr_{i}).
$$
Then,
\[
|\{x \in \mathbb{R}^{n}: |S_{1, \psi}b^{(N)}(x)|\geq\frac{\rho }{2\sqrt{C_{2}}}\}|\le I + II+ III,
\]
where
$$
I:=	  |\Omega\cup E^{*}|,
$$
$$
II:=|\{x \in \mathbb{R}^{n}\setminus(\Omega\cup E^{*}): |S_{1, \psi}(b^{N}-\rho 1_{E})(x)|\geq\frac{\rho }{4\sqrt{C_{2}}}\}|,
$$
and
$$
III:=|\{x \in \mathbb{R}^{n}: |S_{1, \psi} 1_{E}(x)|\geq\frac{1}{4\sqrt{C_{2}}}\}|.
$$

The estimate of $I$ results from Lemma \ref{cube}, properties (2)-(3) and definition of $E$
\[
I\le |\Omega|+|E^{*}|\lesssim |\Omega|+|E| \le \frac{1}{\rho}\|f\|_{L^{1}(\mathbb{R}^{n})}+\frac{1}{\rho}\|b\|_{L^{1}(\mathbb{R}^{n})}\lesssim \frac{1}{\rho}\|f\|_{L^{1}(\mathbb{R}^{n})}.
\]

We know that
$$
\supp(b_{i}-\rho 1_{E_{i}})\subset Q_{i}\cup Q(c_{i},r_{i}), \quad \int_{\mathbb{R}^{n}}(b_{i}(x)-\rho 1_{E_{i}}(x))dx=0,
$$
and
\begin{align*}
Q(c_{i},6n r_{i})\subset 6n Q_{i}\cup Q(c_{i},6n r_{i})\subset \Omega \cup E^{*}.
\end{align*}
Moreover, because of the disjointness of the $Q_{i}'s$ and $E_{i}'s,$ the support of the function $b^{N}-\rho 1_{E}$ is $Q_{j}\cup Q(c_{j},r_{j}).$ By the geometry of the cubes, we have
\begin{equation*}
Q_{j}\cup Q(c_{j},r_{j})=
\begin {cases}
 Q(c_{j},r_{j}) & Q_{j}\subset Q(c_{j},r_{j}) \\ 
 Q_{j}& Q(c_{j},r_{j})\subset Q_{j}.
\end {cases}
\end{equation*} 

It is worth mentioning that $$\int_{\mathbb{R}^{n}}(b_{i}(x)-\rho 1_{E_{i}}(x))dx=0.$$ In fact, by using $|E_{i}|=\frac{\int_{Q_{i}}b_{i}(x)dx}{\rho}$ we see that
$$
\int_{\mathbb{R}^{n}}b_{i}(x)dx-\rho\int_{\mathbb{R}^{n}}1_{E_{i}}(x)dx=\int_{\mathbb{R}^{n}}b_{i}(x)dx-\rho|E_{i}|=\int_{\mathbb{R}^{n}}b_{i}(x)dx-\dfrac{\rho \int_{Q_{i}} b_{i}(x)dx}{\rho}=0.
$$
 Thanks to Chebychev's inequality and disjointness of $Q_{i}'s$ and $E_{i}'s$ as well as Lemma \ref{norm1}, $II$ is estimated as
 \begin{align*}
II&\lesssim \frac{1}{\rho}\int_{\mathbb{R}^{n}\setminus \Omega \cup E^{*}}|S_{1,\psi}(b^{N}-\rho 1_{E})(x)|dx
\\
\\
&\leq \frac{1}{\rho}\int_{\mathbb{R}^{n}\setminus Q(c_{j},6n r_{j})}|S_{1,\psi}(b^{N}-\rho 1_{E})(x)|dx\\
&\le \frac{1}{\rho} \|b^{N}-\rho 1_{E}\|_{L^{1}(\mathbb{R}^{n})}\\
&\leq\frac{1}{\rho} (\|b^{N}\|_{L^{1}(\mathbb{R}^{n})}+\|\rho 1_{E})\|_{L^{1}(\mathbb{R}^{n})})\\
&\leq \frac{1}{\rho} (\|b\|_{L^{1}(\mathbb{R}^{n})}+ \rho|E|)\\
&\leq \frac{1}{\rho}\|f\|_{L^{1}(\mathbb{R}^{n})}+\sum^{N}_{i=1}|E_{i}|\\
&= \frac{1}{\rho}\|f\|_{L^{1}(\mathbb{R}^{n})}+\sum^{N}_{i=1}\frac{\int_{Q_{i}}b_{i}(x)dx}{\rho}\\
&=\frac{1}{\rho}\|f\|_{L^{1}(\mathbb{R}^{n})}+\frac{\int_{Q_{i}}b^{N}(x)dx}{\rho}\\
&\leq \frac{1}{\rho}\|f\|_{L^{1}(\mathbb{R}^{n})}+\frac{\int_{\mathbb{R}^{n}}b(x)dx}{\rho}\\
&\leq \frac{1}{\rho}\|f\|_{L^{1}(\mathbb{R}^{n})}+\frac{1}{\rho}\|f\|_{L^{1}(\mathbb{R}^{n})}\\
&\lesssim \frac{1}{\rho}\|f\|_{L^{1}(\mathbb{R}^{n})}.
 \end{align*}

  Chebshev's inequality, the boundedness of $S_{1,\psi}$ on $L^{2}(\mathbb{R}^{n})$ and the fact that $\sum^{N}_{i=1}|E_{i}|\le \frac{1}{\rho}\|f\|_{L^{1}(\mathbb{R}^{n})}$ help $III$ to be estimated as follows
\begin{align*}
	III&\lesssim \int_{\mathbb{R}^{n}} |S_{1,\psi}(1_{E})(x)| ^{2}dx\\
	&\lesssim \int_{\mathbb{R}^{n}}|1_{E}(x)|^{2}dx\\
	&\leq |E|\\
	&\le \sum^{N}_{i=1}|E_{i}|\\ &\leq\frac{1}{\rho}\|f\|_{L^{1}(\mathbb{R}^{n})}.
\end{align*}
Putting all estimates together, we get
	\begin{equation*}
	|\{x\in \mathbb{R}^{n}:|S_{1,\psi}f(x)|>\rho\}|\lesssim \frac{1}{\rho}\|f\|_{L^{1}(\mathbb{R}^{n})}.
\end{equation*}

Since continuous functions with compact supports are dense in $L^{1}(\mathbb{R}^{n}),$ the proof for all $f \in L^{1}$ is valid.
\end{proof}

\section{$(1,1)-$ Weak type boundedness of $g^{*}_{\lambda,\psi}$ operator }

In this section, we investigate the $(1,1)-$weak type boundedness of $g^{*}_{\lambda,\psi}$ operator using the $(1,1)-$weak type boundedness of $S_{1,\psi}$ operator proved in previous section.
The proof of result in this section is essentially based on proof of \cite[Theorem 1.1]{huang}.

\begin{theorem}
	For all $\lambda>2,$ we have
	\begin{equation}
		\|g^{*}_{\lambda,\psi}f\|_{L^{1,\infty}}\lesssim\|f\|_{L^{1}}.
	\end{equation}
The implicit constant is independent of $f.$
\end{theorem}
\begin{proof}
	We are going to show that $$\sup_{\xi>0}\xi |\{x\in\mathbb{R}^{n}:g^{*}_{\lambda,\psi}f(x)>\xi\}|\lesssim \|f\|_{L^{1}}.$$
Fix $\xi>0.$ Let $\mathcal{O}=\{x\in\mathbb{R}^{n}: S_{1,\psi}f(x)>1\}$ and $\mathcal{U}=\{x\in\mathbb{R}^{n}:M(\chi_{\mathcal{O}})(x)>\frac{\xi}{2}\},$ where $M$ is the Hardy-Littlewood maximal function and $\chi_{\mathcal{O}}$ is the characteristic function of $\mathcal{O}.$
For $k\in\mathbb{N},$ we set
\[
\mathcal{O}_{k}=\{x\in \mathbb{R}^{n}: S_{1,\psi}f(x)>2^{kn}\xi\}
\]
\[
\mathcal{U}_{k}=\{x\in\mathbb{R}^{n}:M(\chi_{\mathcal{O}_{k}})(x)>\frac{1}{2^{kn+1}}\}.
\]
First, we claim that $\mathcal{U}_{k}\subset \mathcal{U}$ for all $k\in \mathbb{N}.$

If $x\in\mathcal{U}_{k},$ then there is a ball $B(y,r)\subset \mathbb{R}^{n}$ containing $x$ such that  \[\dfrac{|B(y,r)\cap \mathcal{O}_{k}|}{|B(y,r)|}>\frac{1}{2^{kn+1}}.\]

Hence,
\[
\frac{1}{|B(y,r)|}\int_{B(y,r)}S_{1,\psi}f(z)dz\geq\frac{1}{|B(y,r)|}\int_{B(y,r)\cap\mathcal{O}_{k}}S_{1,\psi}f(z)dz\geq\frac{1}{|B(y,r)|}(2^{kn}\xi)|B(y,r)\cap\mathcal{O}_{k}|>\frac{\xi}{2},
\]
whence we obtain $x \in U.$

We decompose the set $\{x\in\mathbb{R}^{n}:g^{*}_{\lambda,\psi}f(x)>\xi\}$ as follows
\[\{x\in\mathbb{R}^{n}:g^{*}_{\lambda,\psi}f(x)>\xi\}=\{x\in \mathcal{U}:g^{*}_{\lambda,\psi}f(x)>\xi\}\cup \{x\in \mathbb{R}^{n}\setminus\mathcal{U}:g^{*}_{\lambda,\psi}f(x)>\xi\}.\]
Hence,
\[|\{x\in\mathbb{R}^{n}:g^{*}_{\lambda,\psi}f(x)>\xi\}|=|\{x\in \mathcal{U}:g^{*}_{\lambda,\psi}f(x)>\xi\}|+|\{x\in \mathbb{R}^{n}\setminus\mathcal{U}:g^{*}_{\lambda,\psi}f(x)>\xi\}|\]
\[=:S_{1}+S_{2}.\]
We estimate $S_{1}$ using $(1,1)-$weak type boundedness of $M$ and $S_{1,\alpha}$ respectively.
\begin{align*}S_{1}=|\{x\in \mathcal{U}:g^{*}_{\lambda,\psi}f(x)>\xi\}|\leq |\mathcal{U}|&=|\{x\in\mathbb{R}^{n}:M(\chi_{\mathcal{O}})(x)>\frac{1}{2\alpha^{n}}\}|\\
&\leq\frac{\|\chi_{\mathcal{O}}\|_{L^{1}}}{\frac{1}{2\alpha^{n}}}\\
&=2\alpha^{n}|\mathcal{O}|\\
&=2\alpha^{n}|\{x\in\mathbb{R}^{n}: S_{1,\psi}f(x)>\xi\}|\\
&\leq 2\alpha^{n} \frac{\|f\|_{L^{1}}}{\xi}\\
&\lesssim \frac{\|f\|_{L^{1}(\mathbb{R}^{n})}}{\xi}.
\end{align*}
From here we conclude that
\begin{equation}\label{e1}
	\xi|\{x\in \mathcal{U}:g^{*}_{\lambda,\psi}f(x)>\xi\}|\lesssim \|f\|_{L^{1}(\mathbb{R}^{n})}.
\end{equation}
From \cite{huang}, we recall that for $\lambda>1$
\begin{equation}\label{gs}
g^{*}_{\lambda,\psi}f(x)\lesssim \sum^{\infty}_{k=0}2^{\frac{-k\lambda n}{2}}S_{2^{k},\psi}f(x),~~~x\in\mathbb{R}^{n}.
\end{equation}
Note that by \eqref{gs}, we get
\begin{equation}\label{gs2}
g^{*}_{\lambda,\psi}f^{2}\lesssim \sum^{\infty}_{k=0}2^{-k\lambda n}S_{2^{k},\psi}f^{2}.
\end{equation}

We estimate $S_{2}$ using Chebychev's inequality, relation \eqref{gs2},  Lemma \ref{2},  Cavalieri’s principle and $L^{1,\infty}$ norm  as well as  the $(1.1)-$weak type boundedness of $S_{1,\psi}$ operator as follows
\begin{align*}
S_{2}&=|\{x\in \mathbb{R}^{n}\setminus\mathcal{U}:g^{*}_{\lambda,\psi}f(x)>\xi\}|\\
&=\int_{\{x\in \mathbb{R}^{n}\setminus \mathcal{U}:g^{*}_{\lambda,\psi}f(x)>\xi\}}1dx\\
&\leq \int_{ \mathbb{R}^{n}\setminus \mathcal{U}}\frac{g^{*}_{\lambda,\psi}f(x)^{2}}{\xi^{2}}\\
&=\frac{1}{\xi^{2}} \int_{ \mathbb{R}^{n}\setminus \mathcal{U}}g^{*}_{\lambda,\psi}f(x)^{2}dx\\
&\lesssim \frac{1}{\xi^{2}} \sum^{\infty}_{k=0}2^{-k\lambda n}\int_{ \mathbb{R}^{n}\setminus \mathcal{U}_{k}}S_{2^{k},\psi}f(x) ^{2}dx\\
&\leq \frac{1}{\xi^{2}}\sum^{\infty}_{k=0}2^{-kn(\lambda-1)}\int_{ \mathbb{R}^{n}\setminus \mathcal{O}_{k}}S_{1,\psi}f(x) ^{2}dx\\
&\leq\frac{1}{\xi^{2}}\sum^{\infty}_{k=0}2^{-kn(\lambda-1)}\int^{2^{kn}\xi}_{0}\eta|\{x\in \mathbb{R}^{n}:S_{1,\psi}f(x)>\eta\}|d\eta\\
&\leq \frac{1}{\xi} \sum^{\infty}_{k=0}2^{kn(2-\lambda)} \|S_{1,\psi}f\|_{L^{1,\infty}}\\
&\lesssim  \frac{1}{\xi}\|S_{1,\psi}f\|_{L^{1,\infty}}\\
&\lesssim  \frac{1}{\xi}\|f\|_{L^{1}}.
\end{align*}
Hence,
\begin{equation}\label{e2}
	\xi|\{x\in \mathbb{R}^{n}\setminus\mathcal{U}:g^{*}_{\lambda,\psi}f(x)>\xi\}|\lesssim \|f\|_{L^{1}}.
\end{equation}
Putting all estimates \eqref{e1} and \eqref{e2} together, then we take supremum over all $\xi>0$ from them, it follows that
\[\|g^{*}_{\lambda,\psi}f\|_{L^{1,\infty}}\lesssim \|f\|_{L^{1}}.\]
\end{proof}


\bibliographystyle{amsplain}

\begin{thebibliography}{10}

\bibitem{BuiHormozi}
T. A. Bui, M. Hormozi,
{Weighted bounds for multilinear square functions},
Potential Anal. \textbf{46} (2017), 135--148.

\bibitem{CHI}
Cao, M., Hormozi, M., Iba\"nez-Firnkorn, G. et al. Weak and Strong Type Estimates for the Multilinear Littlewood–Paley Operators. J Fourier Anal Appl 27, 62 (2021).

\bibitem{GM}
A. Ghorbanalizadeh, M. Mikaeili Nia, Strong and Weak type estimate for Littlewood-Paley operators associated with
Laplace-Bessel differential operator, Banach J. Math. Anal., No.  21(2022).
\bibitem{grafakous}
L. Grafakos, Classical Fourier Analysis, 2nd edition, Springer, 2008.
\bibitem{HYY}
M. Hormozi, Y. Sawano, K. Yabuta. "Pointwise domination and weak $ L^ 1$ boundedness of Littlewood-Paley Operators via sparse operators." arXiv preprint arXiv:2111.08241 (2021).
\bibitem{huang}
Y. Huang, Remarks on Weak-Type Estimates for Certain Grand Square Functions. Math Notes 103, 589–592 (2018).
\bibitem{Ler}
A.K. Lerner, {Sharp weighted norm inequalities for Littlewood–Paley operators and singular integrals.} Advances in Mathematics {\bf 226}(5) (2011), 3912-3926.

\bibitem{Ler5}A.K. Lerner, {On sharp aperture-weighted estimates for square functions},
J. Fourier Anal. Appl.
{\bf 20}(4) (2014),  784--800.

\bibitem{Nazarov}
F. Nazarov, S. Treil, A. Volberg, Weak type estimates and Cotlar inequalities for Calder\'{o}n-Zygmund
operators on nonhomogeneous spaces, Internat. Math. Res. Notices {\bf 9} (1998), 463–487.
\bibitem{shi}
Sh. Shi, Q. Xue, K. Yabuta, On the boundedness of multilinear Littlewood--Paley $g^{*}_{\lambda}$ function, Journal de Math\'{e}matiques Pures et Appliqu\'{e}es, Volume 101, Issue 3, March 2014.
 function
\bibitem{cody}
C. B. Stockdale, A different approach to endpoint weak-type estimates for Calder\'{o}n-Zygmund operators, J. Math. Anal. Appl. 487 (2020), no. 2, 124016.
\bibitem{torch}
A. Torchinsky, Real-variable methods in harmonic analysis, in Pure Appl. Math. (Academic Press, Orlando, FL, 1986), Vol. 123.
\end{thebibliography}

\end{document}